\documentclass[12pt]{amsart}
\usepackage{amscd,amssymb,amsthm}
\usepackage{graphicx}

\newtheorem{theorem}{Theorem}[section]
\newtheorem{remark}{Remark}[section]
\newtheorem{example}{Example}[section]
\newtheorem{property}{Property}[section]

\newtheorem{proposition}{Proposition}[section]
\numberwithin{equation}{section}

\title[MEAN DIVERGENCES]{ON MEAN DIVERGENCE MEASURES}

\author{Inder Jeet Taneja}
\address{Departamento de Matem\'{a}tica\\ Universidade Federal
de Santa Catarina\\
88.040-900 Florian\'{o}polis, SC, Brazil}

\email{taneja@mtm.ufsc.br}

\urladdr{http://www.mtm.ufsc.br/$\sim $taneja}

\keywords{Arithmetic Mean; Geometric Mean; Harmonic Mean; Square-
Root Mean; Divergence Measures; Csisz\'{a}r's f-divergence;
Information inequalities.}

\subjclass[2000]{94A17; 26D13}

\begin{document}

\begin{abstract}
\textit{Arithmetic}, \textit{geometric} and \textit{harmonic means}
are the three \textit{classical means} famous in the literature.
Another mean such as \textit{square-root mean} is also known. In
this paper, we have constructed divergence measures based on
nonnegative differences among these means, and established an
interesting inequality by use of properties of \textit{Csisz\'{a}r's
f-divergence}. An improvement over this inequality is also
presented. Comparison of new \textit{mean divergence measures} with
\textit{classical divergence measures} such as \textit{J-divergence}
\cite{jef, kul}, \textit{Jensen-Shannon difference divergence
measure} \cite{bra, sib} and \textit{arithmetic-geometric mean
divergence measure} \cite{tan2} are also established.
\end{abstract}

\maketitle

\section{Generalized Mean of Order $t$}

Let us consider the following well-known mean of order $t$
\begin{equation}
\label{eq1}
M_t (a,b) = \left\{ {{\begin{array}{*{20}c}
 {\left( {\frac{a^t + b^t}{2}} \right)^{1 / t},} & {t \ne 0,} \\
 {\sqrt {ab} ,} & {t = 0,} \\
 {\max \{a,b\},} & {t = \infty ,} \\
 {\min \{a,b\},} & {t = - \infty ,} \\
\end{array} }} \right.
\end{equation}

\noindent for all $a,b \in \mathbb{R}$.

It is also well known (ref. Beckenbach and Bellman \cite{beb}) that
the $M_t (a,b)$ is monotonically non-decreasing function in relation
to $t$. This allow us to conclude the following inequality
\begin{equation}
\label{eq2}
M_{ - \infty } (a,b) \leqslant M_{ - 1} (a,b) \leqslant M_0 (a,b) \leqslant
M_1 (a,b) \leqslant M_2 (a,b) \leqslant M_\infty (a,b),
\end{equation}

\noindent where
\begin{align}
M_{ - 1} (a,b) & = H(a,b) = \frac{2ab}{a + b}  & -
\,\,Harmonic\,\, mean;\notag\\
M_0 (a,b) & = G(a,b) = \sqrt {ab} & - \,\,Geometric \,\,mean;\notag\\
M_1 (a,b) & = A(a,b) = \frac{a + b}{2} & - \,\,Arithmetic\,\,
mean; \notag\\
\intertext{and} M_2 (a,b) & = S(a,b) = \sqrt {\frac{a^2 + b^2}{2}} &
- \,\,Square\,\,root\,\, mean.\notag
\end{align}

In view of this we have the following inequality.
\begin{equation}
\label{eq3}
H(a,b) \leqslant G(a,b) \leqslant A(a,b) \leqslant S(a,b).
\end{equation}

Recently, author \cite{tan8} improved the above inequality
(\ref{eq3}). Also see S\'{a}ndor \cite{san} for different kinds of
inequalities among the means.

Let us consider now the following non-negative differences arising
due to inequality (\ref{eq3}).
\begin{align}
M_{SA} (a,b) & = S(a,b) - A(a,b)
 = \sqrt {\frac{a^2 + b^2}{2}} - \frac{a + b}{2},\notag\\
M_{SG} (a,b) &= S(a,b) - G(a,b)
 = \sqrt {\frac{a^2 + b^2}{2}} - \sqrt {ab},\notag\\
M_{SH} (a,b) & = S(a,b) - H(a,b)
 = \sqrt {\frac{a^2 + b^2}{2}} - \frac{2ab}{a + b},\notag\\
M_{AH} (a,b) & = A(a,b) - H(a,b)
 = \frac{a + b}{2} - \frac{2ab}{a + b},\notag\\
M_{AG} (a,b) &= A(a,b) - G(a,b)
 = \frac{a + b}{2} - \sqrt {ab}\notag
\intertext{and} M_{GH} (a,b) &= G(a,b) - H(a,b)
 = \sqrt {ab} - \frac{2ab}{a + b}.\notag
\end{align}

In view of (\ref{eq2}), we have the following inequalities among
then \textit{mean difference measures}:
\begin{equation}
\label{eq4}
0 \leqslant M_{SA} (a,b) \leqslant M_{SG} (a,b) \leqslant M_{SH} (a,b)
\end{equation}

\noindent and
\begin{equation}
\label{eq5}
0 \leqslant M_{AG} (a,b) \leqslant M_{AH} (a,b).
\end{equation}

\section{Mean Difference Divergence Measures}

Let
\[
\Gamma _n = \left\{ {P = (p_1 ,p_2 ,...,p_n )\left| {p_i >
0,\sum\limits_{i = 1}^n {p_i = 1} } \right.} \right\}, \,\, n
\geqslant 2,
\]

\noindent be the set of all complete finite discrete probability
distributions.\\

Let us take $a = p_i $ and $b = q_i $ in the differences given above
and sum over all $i = 1,2,..,n$, then for all $P,Q \in \Gamma _n $,
we have the following \textit{mean divergence measures}:\\

\textbf{$\bullet$ Square root - arithmetic mean divergence}
\[
M_{SA} (P\vert \vert Q) = \sum\limits_{i = 1}^n {\sqrt {\frac{p_i^2 + q_i^2
}{2}} } - 1.
\]

\bigskip
\textbf{$\bullet$ Square root - geometric mean divergence}
\[
M_{SG} (P\vert \vert Q) = \sum\limits_{i = 1}^n {\left( {\sqrt {\frac{p_i^2
+ q_i^2 }{2}} - \sqrt {p_i q_i } } \right)} .
\]

\bigskip
\textbf{$\bullet$ Square root - harmonic mean divergence}
\[
M_{SH} (P\vert \vert Q) = \sum\limits_{i = 1}^n {\left( {\sqrt {\frac{p_i^2
+ q_i^2 }{2}} - \frac{2p_i q_i }{p_i + q_i }} \right)} .
\]

\bigskip
\textbf{$\bullet$ Arithmetic -- geometric mean divergence}
\[
M_{AG} (P\vert \vert Q) = 1 - \sum\limits_{i = 1}^n {\sqrt {p_i q_i } } .
\]

\bigskip
\textbf{$\bullet$ Arithmetic -- harmonic mean divergence}
\[
M_{AH} (P\vert \vert Q) = 1 - \sum\limits_{i = 1}^n {\frac{2p_i q_i }{p_i +
q_i }} .
\]

\bigskip
\textbf{$\bullet$ Geometric -- harmonic mean divergence}
\[
M_{GH} (P\vert \vert Q) = \sum\limits_{i = 1}^n {\left( {\sqrt {p_i q_i } -
\frac{2p_i q_i }{p_i + q_i }} \right)} = \sum\limits_{i = 1}^n {\frac{\sqrt
{p_i q_i } \left( {\sqrt {p_i } - \sqrt {q_i } } \right)^2}{p_i + q_i }} .
\]

After simplification, we can write
\[
M_{AG} (P\vert \vert Q) = 1 - B(P\vert \vert Q) = h(P\vert \vert Q) =
\frac{1}{2}\sum\limits_{i = 1}^n {\left( {\sqrt {p_i } - \sqrt {q_i } }
\right)^2} .
\]

\noindent where $B(P\vert \vert Q)$ is the Bhattacharyya \cite{bha}
\textit{ distance }and $h(P\vert \vert Q)$ is the well known
Hellinger \cite{hel} \textit{discrimination}.

Also we can write
\[
M_{AH} (P\vert \vert Q) = 1 - W(P\vert \vert Q) = \frac{1}{2}\Delta (P\vert
\vert Q) = \sum\limits_{i = 1}^n {\frac{(p_i - q_i )^2}{2(p_i + q_i )}} ,
\]

\noindent where $W(P\vert \vert Q)$ is the \textit{harmonic mean
divergence} and $\Delta (P\vert \vert Q)$ is the well known
\textit{triangular discrimination}.

Some studies on \textit{square root - arithmetic mean divergence}
can be seen in \"{O}sterreicher
and Vajda \cite{osv} and Dragomir et al. \cite{dgp}.\\

In view of (\ref{eq4}) and (\ref{eq5}), we have the following
inequalities
\begin{equation}
\label{eq6} 0 \leqslant M_{SA} (P\vert \vert Q) \leqslant M_{SG}
(P\vert \vert Q) \leqslant M_{SH} (P\vert \vert Q) \,\,
\end{equation}

\noindent and
\begin{equation}
\label{eq7}
0 \leqslant h(P\vert \vert Q) \leqslant \frac{1}{2}\Delta (P\vert \vert Q).
\end{equation}

In this paper our aim is to obtain an inequality relating the
\textit{mean divergence measures} given above. This shall be done by
use of \textit{Csisz\'{a}r's }$f - $\textit{divergence}.

\section{Csisz\'{a}r's $f-$Divergence and Mean Divergence Measures}

Given a convex function $f:(0,\infty ) \to \mathbb{R}$, the $f -
$\textit{divergence} measure introduced by Csisz\'{a}r \cite{csi1}
is given by
\begin{equation}
\label{eq8}
C_f (P\vert \vert Q) =
\sum\limits_{i = 1}^n {q_i f\left( {\frac{p_i }{q_i }} \right)} ,
\end{equation}

\noindent for all $P,Q \in \Gamma _n $.

The following theorem is well known in the literature \cite{csi1, csi2}:\\

\begin{property} \label{prt31} Let the function $f:[0,\infty ) \to
\mathbb{R}$ be differentiable convex and normalized, i.e., $f(1) =
0$, then the Csisz\'{a}r $f - $\textit{divergence}, $C_f (P\vert
\vert Q)$ is \textit{nonnegative} and \textit{convex} in the pair of
probability distribution $(P,Q) \in \Gamma _n \times \Gamma _n $.
\end{property}

The \textit{mean divergence measures} given in Section 2 can be
written as examples of (\ref{eq8}) and applying property \ref{prt31}
we can check the \textit{nonnegativity} and \textit{convexity} of
these measures. Here below we shall give these as examples.

\begin{example} \label{exa31} Let us consider
\[
f_{SA} (x) = \sqrt {\frac{x^2 + 1}{2}} - \frac{x + 1}{2}, \,\, x \in
(0,\infty ),
\]

\noindent in (\ref{eq8}), then we have $C_f (P\vert \vert Q) =
M_{SA} (P\vert \vert Q)$.

Moreover,
\[
{f}'_{SA} (x) = \frac{x}{\sqrt 2 \sqrt {x^2 + 1} } - \frac{1}{2},
\]

\noindent and
\[
{f}''_{SA} (x) = \frac{1}{\sqrt 2 (x^2 + 1)\sqrt {x^2 + 1} }.
\]

Thus we have ${f}''_{SA} (x) > 0$ for all $x \in (0,\infty )$. Also,
we have $f_{SA} (1) = 0$. In view of this we can say that
\textit{the square root -- geometric mean divergence} is
\textit{nonnegative} and \textit{convex} in the pair of probability
distributions $(P,Q) \in \Gamma _n \times \Gamma _n $.
\end{example}

\begin{example} \label{exa32} Let us consider
\[
f_{SG} (x) = \sqrt {\frac{x^2 + 1}{2}} - \sqrt x , \,\, x \in
(0,\infty ),
\]

\noindent in (\ref{eq8}), then we have $C_f (P\vert \vert Q) =
M_{SG} (P\vert \vert Q)$.

Moreover,
\[
{f}'_{SG} (x) = \frac{1}{\sqrt 2 }\left( {\frac{x}{\sqrt {x^2 + 1} } -
\frac{1}{\sqrt x }} \right),
\]

\noindent and
\[
{f}''_{SG} (x) = \frac{1}{\sqrt 2 (x^2 + 1)\sqrt {x^2 + 1} } +
\frac{1}{4x\sqrt x }.
\]

Thus we have ${f}''_{SG} (x) > 0$ for all $x \in (0,\infty )$. Also,
we have $f_{SG} (1) = 0$. In view of this we can say that
\textit{the square root -- geometric mean divergence} is
\textit{nonnegative} and \textit{convex} in the pair of probability
distributions $(P,Q) \in \Gamma _n \times \Gamma _n $.
\end{example}

\begin{example} \label{exa33} Let us consider
\[
f_{SH} (x) = \sqrt {\frac{x^2 + 1}{2}} - \frac{2x}{x + 1}, \,\, x
\in (0,\infty ) \,\, ,
\]

\noindent in (\ref{eq8}), then we have $C_f (P\vert \vert Q) =
M_{SH} (P\vert \vert Q)$.

Moreover,
\[
{f}'_{SH} (x) = \frac{x}{\sqrt 2 \sqrt {x^2 + 1} } - \frac{2}{(x + 1)^2},
\]

\noindent and
\[
{f}''_{SH} (x) = \frac{1}{\sqrt 2 (x + 1)^2\sqrt {x^2 + 1} } + \frac{4}{(x +
1)^3}.
\]

Thus we have ${f}''_{SH} (x) > 0$ for all $x \in (0,\infty )$. Also,
we have $f_{SH} (1) = 0$. In view of this we can say that
\textit{the square root -- geometric mean divergence} is
\textit{nonnegative} and \textit{convex} in the pair of probability
distributions $(P,Q) \in \Gamma _n \times \Gamma _n $.
\end{example}

\begin{example} \label{exa34} Let us consider
\[
f_h (x) = \frac{1}{2}(\sqrt x - 1)^2, \,\, x \in (0,\infty ),
\]

\noindent in (\ref{eq8}), then we have $C_f (P\vert \vert Q) =
h(P\vert \vert Q)$.

Moreover,
\[
{f}'_h (x) = \frac{\sqrt x - 1}{2\sqrt x },
\]

\noindent and
\[
{f}''_h (x) = \frac{1}{4x\sqrt x }.
\]

Thus we have ${f}''_h (x) > 0$ for all $x \in (0,\infty )$. Also, we
have $f_h (1) = 0$. In view of this we can say that \textit{the
square root -- geometric mean divergence} is \textit{nonnegative}
and \textit{convex} in the pair of probability distributions $(P,Q)
\in \Gamma _n \times \Gamma _n $.
\end{example}

\begin{example} \label{exa35} Let us consider
\[
f_\Delta (x) = \frac{(x - 1)^2}{x + 1}, \,\, x \in (0,\infty ),
\]

\noindent in (\ref{eq6}), then we have $C_f (P\vert \vert Q) =
\Delta (P\vert \vert Q)$.

Moreover,
\[
{f}'_\Delta (x) = \frac{(x - 1)(x + 3)}{(x + 1)^2},
\]

\noindent and
\[
{f}''_\Delta (x) = \frac{8}{(x + 1)^3}.
\]

Thus we have ${f}''_\Delta (x) > 0$ for all $x \in (0,\infty )$.
Also, we have $f_\Delta (1) = 0$. In view of this we can say that
\textit{the square root -- geometric mean divergence} is
\textit{nonnegative} and \textit{convex} in the pair of probability
distributions $(P,Q) \in \Gamma _n \times \Gamma _n $.
\end{example}

\section{Bounds on Mean Divergence Measures}

In this section we shall give bounds on the measures given in
Section 2. In order to get these bounds we shall make use of the
properties of \textit{Csisz\'{a}r's f-divergence} due to Dragomir
\cite{dra}.

\begin{property} \label{prt41} Let $f:\mathbb{R}_ + \to \mathbb{R}$ be
differentiable convex and normalized i.e., $f(1) = 0$. If $P,Q \in
\Gamma _n $, then we have
\begin{equation}
\label{eq9}
0 \leqslant C_f (P\vert \vert Q) \leqslant E_{C_f } (P\vert \vert Q),
\end{equation}

\noindent where
\[
E_{C_f } (P\vert \vert Q) = \sum\limits_{i = 1}^n {(p_i - q_i )}
{f}'(\frac{p_i }{q_i }).
\]
\end{property}

\subsection*{4.1. Square root -- arithmetic mean divergence} We have
\[
0 \leqslant M_{SA} (P\vert \vert Q) \leqslant E_{SA} (P\vert \vert Q),
\]

\noindent where
\begin{align}
E_{SA} (P\vert \vert Q) & =
\sum\limits_{i = 1}^n {(p_i - q_i )\left( {\frac{p_i }{\sqrt
{2(p_i^2 + q_i^2 )} }} \right)} \notag\\
& = \sum\limits_{i = 1}^n {\frac{p_i^2 - p_i q_i + q_i^2 - q_i^2
}{\sqrt {2(p_i^2 + q_i^2 )} }}\notag\\
& = M_{SA} (P\vert \vert Q) + 1 - \sum\limits_{i = 1}^n {\frac{q_i
(p_i + q_i )}{\sqrt {2(p_i^2 + q_i^2 )} }}\notag\\
& = M_{SA} (P\vert \vert Q) + \sum\limits_{i = 1}^n {q_i \left[
{\frac{\sqrt {2(p_i^2 + q_i^2 )} - (p_i + q_i )}{\sqrt {2(p_i^2 +
q_i^2 )} }} \right]}\notag\\
& = M_{SA} (P\vert \vert Q) + \xi _{SA} (P\vert \vert Q), \notag
\end{align}
\noindent with
\[
\xi _{SA} (P\vert \vert Q) = \sum\limits_{i = 1}^n {\frac{\sqrt 2 q_i
}{\sqrt {p_i^2 + q_i^2 } }\left[ {\sqrt {\frac{p_i^2 + q_i^2 }{2}} -
\frac{p_i + q_i }{2}} \right]} .
\]

In view of (\ref{eq9}), we can say that $\xi _{SA} (P\vert \vert Q)
\geqslant 0$.

\subsection*{4.2. Square root -- geometric mean divergence} We have
\[
0 \leqslant M_{SG} (P\vert \vert Q) \leqslant E_{SG} (P\vert \vert Q),
\]

\noindent where
\begin{align}
E_{SG} (P\vert \vert Q) &= \sum\limits_{i = 1}^n {(p_i - q_i )\left(
{\frac{p_i }{\sqrt {2(p_i^2 + q_i^2 )} } - \frac{\sqrt {q_i }
}{2\sqrt {p_i } }} \right)}\notag\\
& = \sum\limits_{i = 1}^n {\left( {\frac{p_i^2 - p_i q_i + q_i^2 -
q_i^2 }{\sqrt {2(p_i^2 + q_i^2 )} } - \frac{\sqrt {q_i (p_i - q_i )}
}{2\sqrt {p_i } }} \right)}\notag\\
& = M_{SG} (P\vert \vert Q) + \sum\limits_{i = 1}^n {\left(
{\frac{\sqrt {q_i } (p_i + q_i )}{\sqrt {p_i } }\frac{q_i (p_i + q_i
)}{\sqrt {2(p_i^2 + q_i^2 )} }} \right)}\notag\\
& = M_{SG} (P\vert \vert Q) + \sum\limits_{i = 1}^n {\sqrt {q_i }
(p_i + q_i )\left( {\frac{\sqrt {p_i^2 + q_i^2 } - \sqrt {2p_i q_i }
}{2\sqrt {p_i } \sqrt {p_i^2 + q_i^2 } }} \right)}\notag\\
& = M_{SG} (P\vert \vert Q) + \sum\limits_{i = 1}^n {\frac{\sqrt
{q_i } (p_i + q_i )}{2\sqrt {p_i } \sqrt {p_i^2 + q_i^2 } }\left(
{\sqrt {p_i^2 + q_i^2 } - \sqrt {2p_i q_i } } \right)}\notag \\
& = M_{SG} (P\vert \vert Q) + \xi _{SG} (P\vert \vert Q),\notag
\end{align}
\noindent with
\[
\xi _{SG} (P\vert \vert Q) = \sum\limits_{i = 1}^n {(p_i + q_i )\sqrt
{\frac{q_i }{2p_i (p_i^2 + q_i^2 )}} \left( {\sqrt {\frac{p_i^2 + q_i^2
}{2}} - \sqrt {p_i q_i } } \right)} .
\]

In view of (\ref{eq9}), we can say that $\xi _{SG} (P\vert \vert Q)
\geqslant 0$.

\subsection*{4.3. Square root -- harmonic mean divergence} We have
\[
0 \leqslant M_{SH} (P\vert \vert Q) \leqslant E_{SH} (P\vert \vert Q),
\]

\noindent where
\begin{align}
E_{SH} (P\vert \vert Q) &= \sum\limits_{i = 1}^n {(p_i - q_i )\left(
{\frac{p_i }{\sqrt {2(p_i^2 + q_i^2 )} } - \frac{2q_i^2 }{(p_i + q_i
)^2}} \right)}\notag\\
& = \sum\limits_{i = 1}^n {\left( {\frac{p_i^2 - p_i q_i + q_i^2 -
q_i^2 }{\sqrt {2(p_i^2 + q_i^2 )} } - \frac{2q_i^2 (p_i - q_i
)}{(p_i^2 + q_i^2 )^2}} \right)}\notag\\
& = M_{SH} (P\vert \vert Q) + \sum\limits_{i = 1}^n {\left(
{\frac{2p_i q_i (p_i + q_i ) - 2q_i^2 (p_i - q_i )}{(p_i + q_i )^2}
- \frac{q_i (p_i + q_i )}{\sqrt {2(p_i^2 + q_i^2 )} }}
\right)}\notag\\
& = M_{SH} (P\vert \vert Q) + \sum\limits_{i = 1}^n {q_i \left[
{\frac{2(p_i^2 + q_i^2 )}{(p_i + q_i )^2} - \frac{(p_i + q_i
)}{\sqrt {2(p_i^2 + q_i^2 )} }} \right]}\notag\\
& = M_{SH} (P\vert \vert Q) + \sum\limits_{i = 1}^n {q_i \left[
{\frac{\left( {\sqrt {2(p_i^2 + q_i^2 )} } \right)^3 - (p_i + q_i
)^3}{(p_i + q_i )^2\sqrt {2(p_i^2 + q_i^2 )} }} \right]}\notag\\
& = M_{SH} (P\vert \vert Q) + \xi _{SH} (P\vert \vert Q),\notag
\end{align}

\noindent with
\[
\xi _{SH} (P\vert \vert Q) = \sum\limits_{i = 1}^n {q_i \left[ {\frac{\left(
{\sqrt {2(p_i^2 + q_i^2 )} } \right)^3 - (p_i + q_i )^3}{(p_i + q_i )^2\sqrt
{2(p_i^2 + q_i^2 )} }} \right]} .
\]

In view of (\ref{eq9}), we can say that $\xi _{SH} (P\vert \vert Q)
\geqslant 0$.

\subsection*{4.4. Hellinger discrimination} We have
\[
0 \leqslant M_h (P\vert \vert Q) \leqslant E_h (P\vert \vert Q),
\]

\noindent where
\begin{align}
E_h (P\vert \vert Q) &= \sum\limits_{i = 1}^n {\frac{(p_i - q_i
)\left( {\sqrt {p_i } - \sqrt {q_i } } \right)}{2\sqrt {p_i }
}}\notag\\
& = \sum\limits_{i = 1}^n {\frac{\left( {\sqrt {p_i } - \sqrt {q_i }
} \right)^2\left( {\sqrt {p_i } + \sqrt {q_i } } \right)}{2\sqrt
{p_i } }}\notag\\
& = M_h (P\vert \vert Q) + \frac{1}{2}\sum\limits_{i = 1}^n {\sqrt
{\frac{q_i }{p_i }} } \left( {\sqrt {p_i } - \sqrt {q_i } }
\right)^2\notag\\
& = M_h (P\vert \vert Q) + \xi _h (P\vert \vert Q),\notag
\end{align}

\noindent with
\[
\xi _h (P\vert \vert Q) = \frac{1}{2}\sum\limits_{i = 1}^n {\sqrt {\frac{q_i
}{p_i }} } \left( {\sqrt {p_i } - \sqrt {q_i } } \right)^2.
\]

Obviously, $\xi _h (P\vert \vert Q) \geqslant 0$.

\subsection*{4.5. Triangular discrimination} We have
\[
0 \leqslant M_\Delta (P\vert \vert Q) \leqslant E_\Delta (P\vert \vert Q),
\]

\noindent where
\begin{align}
E_\Delta (P\vert \vert Q) &= \sum\limits_{i = 1}^n {\frac{(p_i - q_i
)^2(p_i + 3q_i )}{(p_i + q_i )^2}}\notag\\
& = \Delta (P\vert \vert Q) + 2\sum\limits_{i = 1}^n {q_i \left(
{\frac{p_i - q_i }{p_i + q_i }} \right)^2}\notag\\
& = M_\Delta (P\vert \vert Q) + \xi _\Delta (P\vert \vert Q),\notag
\end{align}

\noindent with
\[
\xi _\Delta (P\vert \vert Q) = 2\sum\limits_{i = 1}^n {q_i \left( {\frac{p_i
- q_i }{p_i + q_i }} \right)^2} .
\]

Obviously, $\xi _\Delta (P\vert \vert Q) \geqslant 0$.

\section{Inequalities among Mean Divergence Measures}

In this section we shall obtain inequalities among the measures
given in Section 2.

\begin{property} \label{prt51} Let $f_1 ,f_2 :I \subset \mathbb{R}_ +
\to \mathbb{R}$ be two convex mappings that are normalized, i.e.,
$f_1 (1) = f_2 (1) = 0$ and suppose the assumptions:

(i) $f_1 $ and $f_2 $ are twice differentiable on $(a,b)$;

(ii) there exists the real constants $\alpha ,\beta $ such that $\alpha <
\beta $ and
\begin{equation}
\label{eq10} \alpha \leqslant \frac{f_1 ^{\prime \prime }(x)}{f_2
^{\prime \prime }(x)} \leqslant \beta , \,\, f_2 ^{\prime \prime
}(x) > 0, \,\, \forall x \in (a,b).
\end{equation}

Then,
\begin{equation}
\label{eq11}
\alpha \mbox{ }C_{f_2 } (P\vert \vert Q) \leqslant C_{f_1 } (P\vert \vert Q)
\leqslant \beta \mbox{ }C_{f_2 } (P\vert \vert Q),
\end{equation}

\noindent and
\begin{align}
\label{eq12} \alpha \left[ {E_{f_2 } (P\vert \vert Q) - C_{f_2 }
(P\vert \vert Q)} \right] & \leqslant E_{f_1 } (P\vert \vert Q) -
C_{f_1 } (P\vert \vert Q)\\
& \leqslant \beta \left[ {E_{f_2 } (P\vert \vert Q) - C_{f_2 }
(P\vert \vert Q)} \right]\notag
\end{align}
\end{property}

\begin{proof} Let us consider the functions
\begin{equation}
\label{eq13} p_\alpha (x) = f_1 (x) - \alpha \, f_2 (x)
\end{equation}

\noindent and
\begin{equation}
\label{eq14} p_\beta (x) = \beta \, f_2 (x) - f_1 (x),
\end{equation}

\noindent where $\alpha $ and $\beta $ are as given by (\ref{eq10}).

Since $f_1 (x)$ and $f_2 (x)$ are normalized, i.e., $f_1 (1) = f_2
(1) = 0$, then $p_\alpha (1) = p_\beta (1) = 0$. Also, the functions
$f_1 (x)$ and $f_2 (x)$ are twice differentiable. Then in view of
(\ref{eq10}), we have
\begin{equation}
\label{eq15} {p}''_\alpha (x) = f_1 ^{\prime \prime }(x) - \alpha \,
f_2 ^{\prime \prime }(x)
 = f_2 ^{\prime \prime }(x)\left( {\frac{f_1 ^{\prime \prime }(x)}{f_2
^{\prime \prime }(x)} - \alpha } \right) \geqslant 0
\end{equation}

\noindent and
\begin{equation}
\label{eq16} {p}''_\beta (x) = \beta \, f_2 ^{\prime \prime }(x) -
f_1 ^{\prime \prime }(x)
 = f_2 ^{\prime \prime }(x)\left( {\beta - \frac{f_1 ^{\prime \prime
}(x)}{f_2 ^{\prime \prime }(x)}} \right) \geqslant 0,
\end{equation}

\noindent for all $x \in (a,b)$.

In view of (\ref{eq15}) and (\ref{eq16}), we can say that the
functions $p_\alpha ( \cdot )$ and $p_\beta ( \cdot )$ are convex on
$(a,b)$.

According to Property \ref{prt31}, we have
\begin{equation}
\label{eq17}
C_{p_\alpha } (P\vert \vert Q) = C_{f_1 - \alpha f_2 } (P\vert \vert Q) =
C_{f_1 } (P\vert \vert Q) - \alpha \mbox{ }C_{f_2 } (P\vert \vert Q)
\geqslant 0,
\end{equation}

\noindent and
\begin{equation}
\label{eq18}
C_{q_\beta } (P\vert \vert Q) = C_{\beta f_2 - f_1 } (P\vert \vert Q) =
\beta \mbox{ }C_{f_2 } (P\vert \vert Q) - C_{f_1 } (P\vert \vert Q)
\geqslant 0.
\end{equation}

Combining (\ref{eq17}) and (\ref{eq18}) we have the proof of
(\ref{eq11}).\\

Now, we shall prove the inequalities (\ref{eq12}). We have seen
above that the real mappings $p_\alpha ( \cdot )$ and $p_\beta (
\cdot )$ defined over $\mathbb{R}_ + $ are normalized, twice
differentiable and convex on $(a,b)$. Applying the $r.h.s.$ of the
inequalities (\ref{eq11}), we have
\begin{equation}
\label{eq19}
C_{p_\alpha } (P\vert \vert Q) \leqslant E_{C_{p_\alpha } } (P\vert \vert
Q)
\end{equation}

\noindent and
\begin{equation}
\label{eq20} C_{q_\beta } (P\vert \vert Q) \leqslant E_{C_{p_{_\beta
} } } (P\vert \vert Q).
\end{equation}

Moreover,
\begin{equation}
\label{eq21}
C_{p_\alpha } (P\vert \vert Q) = C_{f_1 } (P\vert \vert Q) - \alpha \mbox{
}C_{f_2 } (P\vert \vert Q)
\end{equation}

\noindent and
\begin{equation}
\label{eq22}
C_{p_{_\beta } } (P\vert \vert Q) = \beta \mbox{ }C_{f_2 } (P\vert \vert Q)
- C_{f_1 } (P\vert \vert Q).
\end{equation}

In view of (\ref{eq19}) and (\ref{eq21}), we have
\begin{align}
C_{f_1 } (P\vert \vert Q) - \alpha \mbox{ }C_{f_2 } (P\vert \vert Q)
& \leqslant E_{C_{f_1 ^\prime - \alpha f_2 ^\prime } } \left(
{P\vert \vert Q} \right)\notag\\
& = E_{f_1 } (P\vert \vert Q) - \alpha E_{f_2 } (P\vert \vert
Q).\notag
\end{align}

This gives,
\[
\alpha \left[ {E_{C_{f_2 } } (P\vert \vert Q) - C_{f_2 } (P\vert \vert Q)}
\right] \leqslant E_{C_{f_1 } } (P\vert \vert Q) - C_{f_1 } (P\vert \vert
Q).
\]

Thus, we have the $l.h.s.$ of the inequalities (\ref{eq12}).

Again in view of (\ref{eq20}) and (\ref{eq22}), we have
\begin{align}
\beta \mbox{ C}_{f_2 } (P\vert \vert Q) - C_{f_1 } (P\vert \vert Q)
& \leqslant E_{C_{\beta f_2 - f_1 } } \left( {P\vert \vert Q}
\right)\notag\\
&= \beta \mbox{ E}_{\mbox{C}_{f_2 } } (P\vert \vert Q) - E_{C_{f_1 }
} (P\vert \vert Q).\notag
\end{align}

After simplifying, we get
\[
E_{f_1 } (P\vert \vert Q) - C_{f_1 } (P\vert \vert Q) \leqslant \beta \left[
{E_{f_2 } (P\vert \vert Q) - C_{f_2 } (P\vert \vert Q)} \right].
\]

Thus we have the $r.h.s.$ of the inequalities (\ref{eq12}). This
completes the proof of the property.
\end{proof}

Now, we shall apply the above theorem for the measures given in
Section 2.

\begin{theorem} \label{the51} The following inequalities among the
mean difference and auxiliary divergences hold:
\begin{equation}
\label{eq23}
M_{SA} (P\vert \vert Q) \leqslant \frac{1}{3}M_{SH} (P\vert \vert Q)
\leqslant \frac{1}{4}\Delta (P\vert \vert Q) \leqslant \frac{1}{2}M_{SG}
(P\vert \vert Q) \leqslant h(P\vert \vert Q),
\end{equation}

\noindent and
\begin{equation}
\label{eq24}
\xi _{SA} (P\vert \vert Q) \leqslant \frac{1}{3}\xi _{SH} (P\vert \vert Q)
\leqslant \frac{1}{4}\xi _\Delta (P\vert \vert Q) \leqslant \frac{1}{2}\xi
_{SG} (P\vert \vert Q) \leqslant \xi _h (P\vert \vert Q).
\end{equation}
\end{theorem}

The proof is based on the following propositions.

\begin{proposition} \label{pro51} The following inequalities hold:
\begin{equation}
\label{eq25}
0 \leqslant M_{SA} (P\vert \vert Q) \leqslant \frac{1}{3}M_{SH} (P\vert
\vert Q),
\end{equation}

\noindent and
\begin{equation}
\label{eq26} 0 \leqslant \xi _{SA} (P\vert \vert Q) \leqslant
\frac{1}{3}\xi _{SH} (P\vert \vert Q) \,\, .
\end{equation}
\end{proposition}

\begin{proof} Let us consider
\[
g_{SA\_SH} (x) = \frac{{f}''_{SA} (x)}{{f}''_{SH} (x)} = \frac{(x +
1)^3}{(x + 1)^3 + 4\sqrt 2 (x^2 + 1)^{3 / 2}}, \,\, x \in (0,\infty
).
\]

This gives
\begin{equation}
\label{eq27}
{g}'_{SA\_SH} (x) = - \frac{24(x - 1)(x^2 + 1)(x + 1)^2}{\sqrt {2(x^2 + 1)}
\left[ {(x + 1)^3 + 4\sqrt 2 (x^2 + 1)^{3 / 2}} \right]^2}
\left\{ {{\begin{array}{*{20}c}
 { \geqslant 0,} & {x \leqslant 1} \\
 { \leqslant 0,} & {x \geqslant 1} \\
\end{array} }} \right..
\end{equation}

In view of (\ref{eq27}) we conclude that the function $g_{SA\_SH}
(x)$ increasing in $x \in (0,1)$ and decreasing in $x \in (1,\infty
)$, and hence
\begin{equation}
\label{eq28}
\beta = \mathop {\sup }\limits_{x \in (0,\infty )} g_{SA\_SH} (x) =
g_{SA\_SH} (1) = \frac{1}{3}.
\end{equation}

Now (\ref{eq28}) together with (\ref{eq11}) and (\ref{eq12}) give
respectively (\ref{eq25}) and (\ref{eq26}).
\end{proof}

\begin{proposition} \label{pro52} The following inequalities hold:
\begin{equation}
\label{eq29}
0 \leqslant M_{SA} (P\vert \vert Q) \leqslant \frac{1}{4}\Delta (P\vert
\vert Q),
\end{equation}

\noindent and
\begin{equation}
\label{eq30}
0 \leqslant \xi _{SA} (P\vert \vert Q) \leqslant \frac{1}{4}\xi _\Delta
(P\vert \vert Q).
\end{equation}
\end{proposition}

\begin{proof} Let us consider
\[
g_{SA\_\Delta } (x) = \frac{{f}''_{SA} (x)}{{f}''_\Delta (x)} =
\frac{(x + 1)^3}{8\sqrt 2 (x^2 + 1)^{3 / 2}}, \,\, x \in (0,\infty
),
\]

This gives
\begin{equation}
\label{eq31}
{g}'_{SA\_\Delta } (x) = - \frac{3(x - 1)(x + 1)^2}{8\sqrt 2 (x^2 + 1)^{5 /
2}}
\left\{ {{\begin{array}{*{20}c}
 { \geqslant 0,} & {x \leqslant 1} \\
 { \leqslant 0,} & {x \geqslant 1} \\
\end{array} }} \right..
\end{equation}

In view of (\ref{eq31}), we conclude that the function
$g_{SA\_\Delta } (x)$ is increasing in $x \in (0,1)$ and decreasing
in $x \in (1,\infty )$, and hence
\begin{equation}
\label{eq32}
M = \mathop {\sup }\limits_{x \in (0,\infty )} g_{SA\_\Delta } (x) =
g_{SA\_\Delta } (1) = \frac{1}{4}.
\end{equation}

Now (\ref{eq32}) together with (\ref{eq11}) and (\ref{eq12}) give
respectively (\ref{eq29}) and (\ref{eq30}).
\end{proof}

\begin{proposition} \label{pro53} The following inequalities hold:
\begin{equation}
\label{eq33} 0 \leqslant \frac{1}{2}\Delta (P\vert \vert Q)
\leqslant M_{SG} (P\vert \vert Q) \,\, ,
\end{equation}

\noindent and
\begin{equation}
\label{eq34}
0 \leqslant \frac{1}{2}\xi _\Delta (P\vert \vert Q) \leqslant \xi _{SG}
(P\vert \vert Q).
\end{equation}
\end{proposition}

\begin{proof} Let us consider
\[
g_{SG\_\Delta } (x) = \frac{{f}''_{SG} (x)}{{f}''_\Delta (x)} =
\frac{(x + 1)^3\left[ {4x^{3 / 2} + \sqrt 2 \mbox{ }(x^2 + 1)^{3 /
2}} \right]}{32\sqrt 2 \mbox{ }(x^2 + 1)^{3 / 2}x^{3 / 2}}, \,\, x
\in (0,\infty ),
\]

This gives
\begin{equation}
\label{eq35}
{g}'_{SG\_\Delta } (x) = \frac{3(x + 1)^4(x - 1)\left[ {\sqrt 2 \left( {x^2
+ 1} \right)^{5 / 2} - 8x^{5 / 2}} \right]}{64\sqrt 2 \left[ {x(x + 1)}
\right]^{5 / 2}}.
\end{equation}

Since $x^2 + 1 \geqslant 2x$, then from (\ref{eq39}), we conclude
that
\begin{equation}
\label{eq36}
{g}'_{SG\_\Delta } (x)\left\{ {{\begin{array}{*{20}c}
 { \geqslant 0,} & {x \geqslant 1} \\
 { \leqslant 0,} & {x \leqslant 1} \\
\end{array} }} \right..
\end{equation}

In view of (\ref{eq36}), we conclude that the function
$g_{SG\_\Delta } (x)$ is decreasing in $x \in (0,1)$ and increasing
in $x \in (1,\infty )$, and hence
\begin{equation}
\label{eq37}
\alpha = \mathop {\inf }\limits_{x \in (0,\infty )} g_{SG\_\Delta } (x) =
\mathop {\min }\limits_{x \in (0,\infty )} g_{SG\_\Delta } (x) =
\frac{1}{2}.
\end{equation}

Now (\ref{eq37}) together with (\ref{eq11}) and (\ref{eq12}) give
respectively (\ref{eq33}) and (\ref{eq34}).
\end{proof}

\begin{proposition} \label{pro54} We have the following bounds:
\begin{equation}
\label{eq38}
0 \leqslant M_{SG} (P\vert \vert Q) \leqslant 2\mbox{ }h(P\vert \vert Q),
\end{equation}

\noindent and
\begin{equation}
\label{eq39} 0 \leqslant \xi _{SG} (P\vert \vert Q) \leqslant
2\,\,\xi _h (P\vert \vert Q).
\end{equation}
\end{proposition}

\begin{proof} Let us consider
\[
g_{SG\_h} (x) = \frac{{f}''_{SG} (x)}{{f}''_h (x)} = \frac{4x^{3 /
2} + \sqrt 2 \mbox{ }(x^2 + 1)^{3 / 2}}{\sqrt 2 \mbox{ }(x^2 + 1)^{3
/ 2}}, \,\, x \in (0,\infty ).
\]

This gives
\begin{equation}
\label{eq40}
{g}'_{SG\_h} (x) = - \frac{6(x - 1)(x + 1)\sqrt x }{\sqrt 2 (x^2 + 1)^{5 /
2}}
\left\{ {{\begin{array}{*{20}c}
 { \geqslant 0,} & {x \leqslant 1} \\
 { \leqslant 0,} & {x \geqslant 1} \\
\end{array} }} \right..
\end{equation}

In view of (\ref{eq40}), we conclude that the function $g_{SG\_h}
(x)$ is increasing in $x \in (0,1)$ and decreasing in $x \in
(1,\infty )$, and hence
\begin{equation}
\label{eq41}
\beta = \mathop {\sup }\limits_{x \in (0,\infty )} g_{SG\_h} (x) = g_{SG\_h}
(1) = 2.
\end{equation}

Now (\ref{eq41}) together with (\ref{eq11}) and (\ref{eq12}) give
respectively (\ref{eq38}) and (\ref{eq39}).
\end{proof}

The inequalities (\ref{eq25}), (\ref{eq29}), (\ref{eq33}) and
(\ref{eq38}) together give (\ref{eq23}) and the inequalities
(\ref{eq26}), (\ref{eq30}), (\ref{eq34}) and (\ref{eq39}) together
give (\ref{eq24}). This completes the proof of the Theorem
\ref{the51}.

\begin{remark}
\begin{itemize}
\item[(i)] The divergence measure arising due to
\textit{geometric--harmonic mean} is not studied here because it is
not convex.\\
\item[(ii)]The \textit{auxiliary measures} $\xi _{( \cdot )} (P\vert \vert
Q)$ can be written in terms of \textit{Csisz\'{a}r f-divergence},
but they are not necessarily convex.
\end{itemize}
\end{remark}

\section{Comparison with Classical Divergence Measures}

In this section, we shall present some classical divergence
measures. The following \textit{Jensen-Shannon divergence measure}
\cite{bra,sib} is already known in the literature:
\begin{equation}
\label{eq42}
I(P\vert \vert Q) = \sum\limits_{i = 1}^n {\left[ {A(p_i \ln p_i ,q_i \ln
q_i ) - A\left( {p_i ,q_i } \right)\ln A\left( {p_i ,q_i } \right)} \right]}
.
\end{equation}

Taneja \cite{tan2} presented the following \textit{arithmetic and
geometric divergence measure} arising due to arithmetic and
geometric means:
\begin{equation}
\label{eq43}
T(P\vert \vert Q) = \sum\limits_{i = 1}^n {A(p_i ,q_i )} \ln \frac{A(p_i
,q_i )}{G(p_i ,q_i )}.
\end{equation}

Adding (\ref{eq42}) and (\ref{eq43}), we get
\begin{equation}
\label{eq44}
I(P\vert \vert Q) + T(P\vert \vert Q) = 4J(P\vert \vert Q),
\end{equation}

\noindent where $J(P\vert \vert Q)$ is the well known
Jefferys-Kullback-Leibler \cite{kul,jef} \textit{J-divergence} given
by
\begin{equation}
\label{eq45}
J(P\vert \vert Q) = \sum\limits_{i = 1}^n {(p_i - q_i )\ln \left( {\frac{p_i
}{q_i }} \right)} .
\end{equation}

For more studies on the measures (\ref{eq43})-(\ref{eq45}) with
their generalizations and some statistical applications refer to
Taneja \cite{tan1,tan3,tan6,tan7}. For new symmetric divergence measure
refer to Kumar and Chhina \cite{kuc}\\

Recently, author \cite{tan4, tan5} proved an inequality among these
divergence measures given by
\begin{equation}
\label{eq46}
\frac{1}{4}\Delta (P\vert \vert Q) \leqslant I(P\vert \vert Q) \leqslant
h(P\vert \vert Q) \leqslant \frac{1}{8}J(P\vert \vert Q) \leqslant T(P\vert
\vert Q).
\end{equation}

Finally, combining the inequalities (\ref{eq23}) and (\ref{eq46}),
we have the following interesting inequalities:
\begin{equation}
\label{eq47} M_{SA} (P\vert \vert Q) \leqslant \frac{1}{3}M_{SH}
(P\vert \vert Q) \leqslant \frac{1}{4}\Delta (P\vert \vert Q)
\leqslant \frac{1}{2}M_{SG} (P\vert \vert Q)
\end{equation}
\[
\leqslant h(P\vert \vert Q) \leqslant \frac{1}{8}J(P\vert \vert Q)
\leqslant T(P\vert \vert Q),\notag
\]

\noindent and
\begin{equation}
\label{eq48} M_{SA} (P\vert \vert Q) \leqslant \frac{1}{3}M_{SH}
(P\vert \vert Q) \leqslant \frac{1}{4}\Delta (P\vert \vert Q)
\leqslant I(P\vert \vert Q)
\end{equation}
\[
 \leqslant h(P\vert \vert Q) \leqslant \frac{1}{8}J(P\vert \vert Q)
\leqslant T(P\vert \vert Q).
\]

From the inequalities (\ref{eq23}), (\ref{eq47}) and (\ref{eq48}),
we observe that we don't have relation among the measures
\textit{SG--divergence} and \textit{I--divergence}. Let us check
this by applying Property \ref{prt51}

Let us consider
\[
f_I (x) = \frac{x}{2}\ln x + \frac{x + 1}{2}\ln \left( {\frac{2}{x +
1}} \right), \,\, x \in (0,\infty )
\]

\noindent in (\ref{eq8}), the one gets $C_f (P\vert \vert Q) =
I(P\vert \vert Q)$.

Moreover,
\[
{f}'_I (x) = \frac{1}{2}\ln \left( {\frac{2x}{x + 1}} \right)
\]

\noindent and
\[
{f}''_I (x) = \frac{1}{2x(x + 1)}.
\]

Again, let us consider
\[
g_{SG \_I} (x) = \frac{{f}''_{SG} (x)}{{f}''_I (x)} = \frac{\left[
{4x^{3 / 2} + \sqrt 2 \mbox{ }(x^2 + 1)^{3 / 2}} \right]x(x +
1)}{2\sqrt 2 \mbox{ }(x^2 + 1)^{3 / 2}x^{3 / 2}}, \,\, x \in
(0,\infty ).
\]

The first order derivative of the function $g_{SG \_I}(x)$ is given
by
\[
{g}'_{SG\_I} (x) = \frac{(x - 1)\sigma (x)}{4\sqrt 2 \left( {x^2 + 1}
\right)^{5 / 2}x^{3 / 2}},
\]

\noindent where
\[
\sigma (x) = \sqrt 2 \left( {x^2 + 1} \right)^{5 / 2} - 8x^{3 / 2}\left(
{x^2 + 3x + 1} \right).
\]

In order to apply Property \ref{prt51} we must prove that $\sigma
(x)$ is either negative or positive for $x \in (0,\infty )$, but
$\sigma (1) = - 32.0$ and $\sigma (4.25) = 13.87$. This implies that
we are unable to apply the Property \ref{prt51}.

Moreover, if we check the generating functions in both the cases,
still the result don't hold. Let us denote, $a(x) = f_{SA} (x)$,
$b(x) = \frac{1}{3}f_{SH} (x)$, $c(x) = \frac{1}{4}f_\Delta (x)$,
$d(x) = \frac{1}{2}f_{SG} (x)$, $e(x) = f_I (x)$ and $f(x) = f_h
(x)$ for all $x \in (0,\infty )$. Then we have the following values
of these two functions:
\begin{table}[htbp]
\begin{tabular}
{|p{41pt}|p{41pt}|p{41pt}|p{51pt}|p{55pt}|p{55pt}|p{55pt}|} \hline
$x$& 0.1& 10& 1000& 3000& 3800&
3900 \\
\hline
$a(x)$&
0.1606&
1.6063&
206.6071&
620.8204&
786.5058&
807.2165 \\
\hline
$b(x)$&
0.1762&
1.7627&
235.0363&
706.4403&
895.0021&
918.5723 \\
\hline
$c(x)$&
0.1840&
1.8409&
249.2509&
749.2503&
949.2502&
974.2502 \\
\hline
$d(x)$&
0.1972&
1.9720&
337.7421&
1033.2741&
1312.6808&
1347.6332 \\
\hline
$e(x)$&
0.2136&
2.1368&
342.9660&
1035.5640&
1312.7047&
1347.3491 \\
\hline
$f(x)$&
0.2337&
2.3377&
468.8772&
1445.7277&
1838.8558&
1888.0500 \\
\hline
\end{tabular}
\label{tab1}
\end{table}

We observe from the above table that the values of $d(x)$ and $e(x)$
changes in the interval $x \in \left[ {3800,3900} \right]$, before
it $d(x)$ is always smaller than $e(x)$.\\

Let check the same thing by considering particular values of the
probability distributions. Let us consider $n = 2$, $p_1 = t$, $q_1
= 1 - t$, $p_2 = 1 - t$ and $q_2 = t$. Then we can write
\begin{align}
a(t) & = M_{SA} (P\vert \vert Q) = 2\sqrt {\frac{t^2 + (1 -
t)^2}{2}} - 1,\notag\\
b(t) & = \frac{1}{3}M_{SH} (P\vert \vert Q) = \frac{2}{3}\sqrt
{\frac{t^2 + (1 - t)^2}{2}} - \frac{4}{3}t(1 - t),\notag\\
c(t) & = \frac{1}{4}\Delta (P\vert \vert Q) = \frac{1}{2}(2t - 1)^2,\notag\\
d(t) & = \frac{1}{2}M_{SG} (P\vert \vert Q) = \sqrt {\frac{t^2 + (1
- t)^2}{2}} - \sqrt {t(1 - t)},\notag\\
e(t) & = I(P\vert \vert Q) = t\ln (2t) + (1 - t)\log (2 - 2t)\notag
\intertext{and} f(t) & = h(P\vert \vert Q) = \left( {\sqrt t - \sqrt
{1 - t} } \right)^2,\notag
\end{align}

\noindent for all $t \in [0,1]$ with the convention that $0\log 0 =
0$.\\

Let us compare the measures for some particular values of $t$.
\begin{table}[htbp]
\begin{tabular}
{|p{44pt}|p{47pt}|p{47pt}|p{47pt}|p{47pt}|p{47pt}|p{53pt}|}
\hline
$t$&
0.0001&
0.001&
0.01&
0.1&
0.2&
0.4 \\
\hline
$a(t)$&
0.4140&
0.4128&
0.4001&
0.2806&
0.1662&
0.01980 \\
\hline
$b(t)$&
0.4712&
0.4696&
0.4535&
0.3068&
0.1754&
0.01993 \\
\hline
$c(t)$&
0.4998&
0.4980&
0.4802&
0.3200&
0.1800&
0.02000 \\
\hline
$d(t)$&
0.6970&
0.6747&
0.6005&
0.3403&
0.1830&
0.02004 \\
\hline
$e(t)$&
0.6921&
0.6852&
0.6371&
0.3680&
0.1927&
0.02013 \\
\hline
$f(t)$&
0.9800&
0.9367&
0.8010&
0.4000&
0.2000&
0.02020 \\
\hline
\end{tabular}
\label{tab2}
\end{table}

Here we have considered only the values of $t \in (0,1 / 2]$, since
for $t \in [1 / 2,1)$ the values are symmetric. Moreover, all values
are zero for $t = \frac{1}{2}$. From the table we observe that for
each $t$ fixed the values of the functions are monotonically
increasing, except for $t = 0.0001$, $d(t)$ is bigger than $e(t).$\\

From the example above we conclude that we are unable to establish
an inequality having nine measures in a sequence combining
(\ref{eq47}) and (\ref{eq48}).

\section{Refinement Inequalities}

Now we shall improve the inequality (\ref{eq23}). In order to do so,
we shall again consider the following \textit{non-negative
differences}:
\[
D_{f_k } (P\vert \vert Q) = \sum\limits_{i = 1}^n {q_i f_k \left(
{\frac{p_i }{q_i }} \right)} , \,\, k = 1,2,...,10,
\]

\noindent where , for all $x \in (0,\infty )$, we have
\begin{align}
f_1 (x) & = f_{AG} (x) - \frac{1}{2}f_{SG} (x),\notag\\
f_2 (x) & = f_{AG} (x) - \frac{1}{2}f_{AH} (x),\notag\\
f_3 (x) & = f_{AG} (x) - \frac{1}{3}f_{SH} (x), \notag\\
f_4 (x) & = f_{AG} (x) - f_{SA} (x), \notag\\
f_5 (x) & = \frac{1}{2}f_{SG} (x) - \frac{1}{2}f_{AH} (x), \notag\\
f_6 (x) & = \frac{1}{2}f_{SG} (x) - \frac{1}{3}f_{SH} (x), \notag\\
f_7 (x) & = \frac{1}{2}f_{SG} (x) - f_{SA} (x), \notag\\
f_8 (x) & = \frac{1}{2}f_{AH} (x) - \frac{1}{3}f_{SH} (x), \notag\\
f_9 (x) & = \frac{1}{2}f_{AH} (x) - f_{SA} (x),\notag
\intertext{and} f_{10} (x) & = \frac{1}{3}f_{SH} (x) - f_{SA}
(x).\notag
\end{align}

We can easily verify that
\begin{equation}
\label{eq49}
f_1 (x) = \frac{1}{2}f_4 (x) = f_7 (x),
\end{equation}

\noindent and
\begin{equation}
\label{eq50}
f_8 (x) = \frac{1}{3}f_9 (x) = \frac{1}{2}f_{10} (x).
\end{equation}

For all $x \in (0,\infty )$, we can write
\begin{align}
f_1 (x) & = \left( {\frac{\sqrt x - 1}{2}} \right)^2 - \frac{\sqrt
{2(x^2 + 1)} - 2\sqrt x }{4}
 = A - \left( {\frac{G + S}{2}} \right),\notag\\
f_2 (x) & = \left( {\frac{\sqrt x - 1}{2}} \right)^2 - \frac{(x -
1)^2}{4(x + 1)}
 = \left( {\frac{A + H}{2}} \right) - G,\notag\\
f_3 (x) & = \left( {\frac{\sqrt x - 1}{2}} \right)^2 - \frac{\sqrt
{2(x^2 + 1)} }{6} + \frac{2x}{2(x + 1)}
= \frac{1}{3}\left[ {3A + H - \left( {S + 3G} \right)} \right],\notag\\
f_5 (x) & = \frac{\sqrt {2(x^2 + 1)} }{2} - \sqrt x - \frac{(x -
1)^2}{2(x + 1)}
= \frac{1}{2}\left[S + H - (A + G)\right],\notag\\
f_6 (x) & = \frac{\sqrt {2(x^2 + 1)} }{12} - \frac{\sqrt x }{2} +
\frac{2x}{3(x + 1)} = \frac{1}{6}\left[ {S + 2H - 3G} \right],\notag
\intertext{and} f_8 (x) & = \frac{(x - 1)^2}{4(x + 1)} - \frac{\sqrt
{2(x^2 + 1)} }{6} + \frac{2x}{3(x + 1)} = \frac{1}{6}\left[ {3A -
\left( {2S + H} \right)} \right],\notag
\end{align}

\noindent where $A = \frac{x + 1}{2}$, $G = \sqrt x $, $H =
\frac{2x}{x + 1}$ and $S = \sqrt {\frac{x^2 + 1}{2}} $ are
respectively \textit{arithmetic}, \textit{geometric},
\textit{harmonic} and \textit{square-root means} between $x$ and 1.

\begin{theorem} \label{the71} The following inequality among the
\textit{new differences} hold:
\begin{equation}
\label{eq51}
D_{f_8 } (P\vert \vert Q) \leqslant \frac{1}{3}D_{f_1 } (P\vert \vert Q)_1
\leqslant \frac{1}{4}D_{f_3 } (P\vert \vert Q) \leqslant \frac{1}{3}D_{f_2 }
(P\vert \vert Q) \leqslant D_{f_6 } (P\vert \vert Q).
\end{equation}
\end{theorem}

\begin{proof} We shall prove each part of the inequality
separately. These inequalities can be proved on similar lines of
theorem \ref{the51}, but here we shall give a simpler proof.\\

\item[(i)] We can write
\begin{align}
f_6 (x) - \frac{1}{3}f_2 (x) & = \frac{1}{6}\left[ {S + 2H - 3G}
\right]
 - \frac{1}{3}\left[ {\left( {\frac{A + H}{2}} \right) - G} \right]\notag\\
& = \frac{1}{6}\left[ {S + H - (A + G)} \right] = \frac{1}{3}f_5 (x)
\geqslant 0, \,\, \forall x \in (0,\infty )\notag
\end{align}

This prove that $\frac{1}{3}f_2 (x) \leqslant f_6 (x)$, $\forall x
\in (0,\infty )$, and consequently, we get
\begin{equation}
\label{eq52}
\frac{1}{3}D_{f_2 } (P\vert \vert Q) \leqslant D_{f_6 } (P\vert \vert Q).
\end{equation}

\item[(ii)] We can write
\begin{align}
\frac{1}{3}f_2 (x) - \frac{1}{4}f_3 (x) & = \frac{1}{3}\left[
{\frac{A + H}{2} - G} \right]
 - \frac{1}{12}\left[ {3A + H - \left( {S + 3G} \right)} \right]\notag\\
& = \frac{1}{12}\left[ {S + H - \left( {A + G} \right)} \right] =
\frac{1}{6}f_5 (x) \geqslant 0, \,\, \forall x \in (0,\infty
).\notag
\end{align}

This prove that $\frac{1}{4}f_3 (x) \leqslant \frac{1}{3}f_2 (x)$,
$\forall x \in (0,\infty )$, and consequently, we get
\begin{equation}
\label{eq53}
\frac{1}{4}D_{f_3 } (P\vert \vert Q) \leqslant \frac{1}{3}D_{f_2 } (P\vert
\vert Q).
\end{equation}

\item[(iii)] We can write
\begin{align}
\frac{1}{4}f_3 (x) - \frac{1}{3}f_1 (x) & = \frac{1}{12}\left[ {3A +
H - \left( {S + 3G} \right)} \right]
 - \frac{1}{3}\left[ {A - \left( {\frac{G + S}{2}} \right)} \right]\notag\\
& = \frac{1}{12}\left[ {S + H - \left( {A + G} \right)} \right] =
\frac{1}{6}f_5 (x) \geqslant 0, \,\, \forall x \in (0,\infty
).\notag
\end{align}

This prove that $\frac{1}{3}f_1 (x) \leqslant \frac{1}{4}f_3 (x)$,
$\forall x \in (0,\infty )$, and consequently, we get
\begin{equation}
\label{eq54}
\frac{1}{3}D_{f_1 } (P\vert \vert Q) \leqslant \frac{1}{4}D_{f_3 } (P\vert
\vert Q).
\end{equation}

\item[(iv)] We can write
\begin{align}
\frac{1}{3}f_1 (x) - f_8 (x) & = \frac{1}{3}\left[ {A - \left(
{\frac{G + S}{2}} \right)} \right]
 - \frac{1}{6}\left[ {3A - \left( {2S + H} \right)} \right]\notag\\
& = \frac{1}{6}\left[ {S + H - \left( {A + G} \right)} \right] =
\frac{1}{3}f_5 (x) \geqslant 0, \,\, \forall x \in (0,\infty
).\notag
\end{align}

This prove that $f_8 (x) \leqslant \frac{1}{3}f_1 (x)$, $\forall x
\in (0,\infty )$, and consequently, we get
\begin{equation}
\label{eq55} D_{f_8 } (P\vert \vert Q) \leqslant \frac{1}{3}D_{f_1 }
(P\vert \vert Q).
\end{equation}

Combining (\ref{eq52})-(\ref{eq55}) we get the required result.
\end{proof}

\begin{remark}
\begin{itemize}
\item[(i)] Simplifying the inequalities given (\ref{eq51}) and using the
nonnegativity of the expression $f_5(x)$, $\forall x \in (0,\infty
)$, we get the following improvement over the inequalities
(\ref{eq23}):
\begin{align}
\label{eq56} M_{GH} & (P\vert \vert Q) \leqslant M_{SA}(P\vert \vert
Q) \leqslant \frac{1}{3}M_{SH}
(P\vert \vert Q) \leqslant \frac{1}{4}\Delta (P\vert \vert Q)\\
& \leqslant \frac{3\Delta (P\vert \vert Q) + 2M_{SG} (P\vert \vert
Q)}{16} \leqslant \frac{h(P\vert \vert Q) + 3M_{SA} (P\vert \vert
Q)}{4}\notag\\
& \qquad \leqslant \frac{h(P\vert \vert Q) + M_{SH} (P\vert \vert
Q)}{4} \leqslant \frac{6M_{SG} (P\vert \vert Q) + \Delta (P\vert
\vert Q)}{4}\notag\\
& \qquad \qquad \leqslant \frac{1}{2}M_{SG} (P\vert \vert Q)
\leqslant h(P\vert \vert Q) \leqslant \frac{1}{2}\Delta (P\vert
\vert Q)\notag
\end{align}

We observe that the measure $M_{GH} (P\vert \vert Q)$ in not convex
in the pair of probability distributions, but even then we are able
to relate it in the above inequalities.

\item[(ii)] Recently, author \cite{tan5} also gave an improvement over the
inequality (\ref{eq46}):
\begin{align}
\label{eq57} \frac{1}{4}\Delta & (P\vert \vert Q) \leqslant I(P\vert
\vert Q) \leqslant \frac{2}{3}h(P\vert \vert Q) + \frac{1}{12}\Delta
(P\vert \vert Q) \leqslant h(P\vert \vert Q)\\
& \leqslant \frac{1}{16}J(P\vert \vert Q) + \frac{1}{2}I(P\vert
\vert Q) \leqslant \frac{1}{3}T(P\vert \vert Q) +
\frac{2}{3}h(P\vert \vert Q)\notag\\
& \qquad \leqslant \frac{1}{8}J(P\vert \vert Q) \leqslant
\frac{2}{3}T(P\vert \vert Q) + \frac{1}{12}\Delta (P\vert \vert Q)
\leqslant T(P\vert \vert Q).\notag
\end{align}
\end{itemize}
The above inequality also improves the one studied by Dragomir et
al. \cite{dsb}.
\end{remark}


\begin{thebibliography}{99}
\setlength{\itemsep}{5pt}

\bibitem{beb} E.F. BECKENBACH and R. BELLMAN, \textit{Inequalities},
Springer-Verlag, New York, 1971.

\bibitem{bha} A.
BHATTACHARYYA, Some Analogues to the Amount of Information and Their
uses in Statistical Estimation, \textit{Sankhya}, \textbf{8}(1946),
1-14.

\bibitem{bra} J. BURBEA and C.R. RAO, On the Convexity of Some Divergence Measures Based on
Entropy Functions, \textit{IEEE Trans. on Inform. Theory},
\textbf{IT-28}(1982), 489-495.

\bibitem{csi1} I. CSISZ\'{A}R,
Information Type Measures of Differences of Probability Distribution
and Indirect Observations, \textit{Studia Math. Hungarica},
\textbf{2}(1967), 299-318.

\bibitem{csi2} I. CSISZ\'{A}R,
On Topological Properties of $f-$Divergences, \textit{Studia Math.
Hungarica}, \textbf{2}(1967), 329-339.

\bibitem{dra} S. S. DRAGOMIR, Some Inequalities for
the Csisz\'{a}r $\Phi $-Divergence - \textit{Inequalities for
Csisz\'{a}r f-Divergence in Information Theory -} Monograph --
Chapter I -- Article 1 --
\textit{http://rgmia.vu.edu.au/monographs/csiszar.htm.}

\bibitem{dsb} S. S. DRAGOMIR, J. SUNDE and C. BUSE, New inequalities for
jeffreys divergence measure, \textit{Tamsui Oxford Journal of
Mathematical Sciences}, \textbf{16}(2)(2000), 295-309.

\bibitem{dgp} S. S. DRAGOMIR, V. GLUSCEVIC and C.E.M. PEARCE, New
Approximations for $f-$Divergence via Trapezoid and Midpoint
Inequalities \textit{http://rgmia.vu.edu.au, RGMIA Research Report
Collection}, \textbf{5}(4)(2002), Article 14.

\bibitem{hel} E. HELLINGER, Neue Begr\"{u}ndung der Theorie der
quadratischen Formen von unendlichen vielen Ver\"{a}nderlichen,
\textit{J. Reine Aug. Math.}, \textbf{136}(1909), 210-271.

\bibitem{jef} H. JEFFREYS, An Invariant Form for the Prior
Probability in Estimation Problems, \textit{Proc. Roy. Soc. Lon.,
Ser. A}, \textbf{186}(1946), 453-461.

\bibitem{kul} S. KULLBACK
and R.A. LEIBLER, On Information and Sufficiency, \textit{Ann. Math.
Statist.}, \textbf{22}(1951), 79-86.

\bibitem{kuc} P. KUMAR and S.A. CHHINA, A Symmetric Information Divergence
Measure of Csiszár's $f-$Divergence Class and its Bounds,
\textit{Computer and Mathematics with Applications},
\textbf{49}(4)(2005), 575-588.

\bibitem{osv} F. \"{O}STERREICHER and I. VAJDA, A New Class of
Metric Divergences on Probability Spaces and its Applicability in
Statistics, \textit{Ann. Inst. Statist. Math.},
\textbf{55}(3)(2003), 639-653.

\bibitem{san} J. S\'{A}NDOR, On Certain Inequalities for Means - II, J. Math.
Analy. and Appl., 199(1996), 629-634.

\bibitem{sib} R. SIBSON, Information
Radius, \textit{Z. Wahrs. und verw Geb.}, \textbf{(14)}(1969),
149-160.

\bibitem{tan1} I.J. TANEJA, On Generalized Information Measures
and Their Applications, Chapter in: \textit{Advances in Electronics
and Electron Physics}, Ed. P.W. Hawkes, Academic Press,
\textbf{76}(1989), 327-413.

\bibitem{tan2} I.J. TANEJA, New Developments in Generalized Information
Measures, Chapter in: \textit{Advances in Imaging and Electron
Physics}, Ed. P.W. Hawkes, \textbf{91}(1995), 37-136.

\bibitem{tan3} I.J. TANEJA, \textit{Generalized Information Measures and their
Applications} - On-line book: \textit{http://www.mtm.ufsc.br/$\sim
$taneja/book/book.html}, 2001.

\bibitem{tan4} I.J. TANEJA, Generalized Symmetric Divergence Measures and
Inequalities - \textit{RGMIA  Research Report Collection},
http://rgmia.vu.edu.au, \textbf{7}(4)(2004), Art. 9. Available
on-line at: arXiv:math.ST/0501301 v1 19 Jan 2005.

\bibitem{tan5} I.J. TANEJA, Refinement Inequalities Among Symmetric
Divergence Measures – \textit{The Australian Journal of Mathematical
Analysis and Applications}, \textbf{2}(1)(2005), Art. 8, pp. 1-23.
Available on-line at: arXiv:math.PR/0501303 v1 19 Jan 2005.

\bibitem{tan6} I.J. TANEJA, Generalized Arithmetic and Geometric Mean
Divergence Measures and Their Statistical Aspects -- To appear in
\textit{Journal of Interdisciplinary Mathematics}. Available on-line
at: arXiv:math.PR/0501297 v1 19 Jan 2005.

\bibitem{tan7} I.J. TANEJA, Bounds on Triangular Discrimination, Harmonic
Mean and Symmetric Chi-Square Divergences - To appear in
\textit{Journal of Concrete and Applicable Mathematics}. Available
on-line at: arXiv:math.PR/0505238 v1 12 May 2005.

\bibitem{tan8} I.J. TANEJA, Refinement Inequalities Among Means - Available
on-line at: arXiv:math. GM/0505192 v1 10 May 2005.

\end{thebibliography}
\end{document}